\documentclass{llncs}
\usepackage[utf8]{inputenc}
\usepackage{graphicx}
\usepackage[draft]{minted}
\usepackage{amssymb, amsmath, bm}
\usepackage[textsize=small]{todonotes}

\usepackage{float}
\floatstyle{ruled}
\newfloat{code}{thp}{lop}
\floatname{code}{Code Block}

\usepackage{hyperref} 
\hypersetup{
colorlinks=false, 
breaklinks=true,
urlcolor= blue, 
linkcolor= blue,
citecolor=red, 
}

\usepackage{lscape}

\newcommand{\bonmin}{{\sc Bonmin}}
\newcommand{\knitro}{{\sc Knitro}}
\newcommand{\scip}{{\sc Scip}}
\newcommand{\couenne}{{\sc Couenne}}
\newcommand{\juniper}{{\sc Juniper}}
\newcommand{\ipopt}{{\sc Ipopt}}
\newcommand{\cbc}{{\sc Cbc}}
\newcommand{\glpk}{{\sc Glpk}}
\newcommand{\gurobi}{{\sc Gurobi}}
\newcommand{\minotaur}{{\sc Minotaur}}

\title{Juniper: An Open-Source Nonlinear \\ Branch-and-Bound Solver in Julia}

\author{Ole Kröger, Carleton Coffrin, Hassan Hijazi, Harsha Nagarajan}
\institute{Los Alamos National Laboratory, Los Alamos, New Mexico, USA}

\begin{document}

\maketitle

\begin{abstract}
\vspace{-0.5cm}
Nonconvex mixed-integer nonlinear programs (MINLPs) represent a challenging class of optimization problems that often arise in engineering and scientific applications.  Because of nonconvexities, these programs are typically solved with global optimization algorithms, which have limited scalability.  However, nonlinear branch-and-bound has recently been shown to be an effective heuristic for quickly finding high-quality solutions to large-scale nonconvex MINLPs, such as those arising in infrastructure network optimization.  This work proposes \juniper{}, a Julia-based open-source solver for nonlinear branch-and-bound.  Leveraging the high-level Julia programming language makes it easy to modify \juniper{}'s algorithm and explore extensions, such as branching heuristics, feasibility pumps, and parallelization.  Detailed numerical experiments demonstrate that the initial release of \juniper{} is comparable with other nonlinear branch-and-bound solvers, such as \bonmin{}, \minotaur{}, and \knitro{}, illustrating that \juniper{} provides a strong foundation for further exploration in utilizing nonlinear branch-and-bound algorithms as heuristics for nonconvex MINLPs.
\end{abstract}

\section{Introduction}
\label{sec:intro}

Many of the optimization problems arising in engineering and scientific disciplines combine both nonlinear equations and discrete decision variables.  Notable examples include the blending/pooling problem \cite{Audet:2004,MINLP168} and the design and operation of power networks \cite{6407493,7038446,h_minlp_opf} and natural gas networks \cite{ijoc.2016.0697}.
All of these problems fall into the class of mixed-integer nonlinear programs (MINLPs), namely,
\begin{equation}
    \tag{MINLP} \label{eq:minlp}
\begin{aligned}
    & \mbox{minimize: } f(x,y) \\
    & \mbox{s.t.} \nonumber \\
    & g_c(x,y) \leq 0 \;\; \forall c \in {\cal C} \\
    & x \in \mathbb{R}^m, y \in \mathbb{Z}^n
\end{aligned}
\end{equation}
where $f$ and $g$ are twice continuously differentiable functions and $x$ and $y$ represent real and discrete valued decision variables, respectively \cite{MINLP2013}.  Combining nonlinear functions with discrete decision variables makes \ref{eq:minlp}s a broad and challenging class of mathematical programs to solve in practice.  To address this challenge, algorithms have been designed for special subclasses of \ref{eq:minlp}s, such as when $f$ and $g$ are convex functions \cite{bonmin,pajarito} or when $f$ and $g$ are nonconvex quadratic functions \cite{scip,boxqp_archive}.  For generic nonconvex functions, global optimization algorithms \cite{brach_and_reduce,couenne,pod1,pod2} are required to solve \ref{eq:minlp}s with a proof of optimality.  However, the scalability of such algorithms is limited and remains an active area of research.  Although global optimization algorithms have been widely successful at solving industrial \ref{eq:minlp}s with a few hundred variables, their limited scalability precludes application to larger real-world problems featuring thousands of variables and constraints, such as AC optimal transmission switching \cite{6939776}.

One approach to addressing the challenge of solving large-scale industrial \ref{eq:minlp}s is to develop heuristics that attempt to quickly find high-quality feasible solutions without guarantees of global optimality. To that end, it has been recently observed that nonlinear branch-and-bound (NLBB) algorithms can be effective heuristics for the nonconvex \ref{eq:minlp}s arising in infrastructure systems \cite{7038446,h_minlp_opf,ijoc.2016.0697} and that they present a promising avenue for solving such problems on real-world scales. To the best of our knowledge, \bonmin{} and \minotaur{}  are the only open-source solvers that implement NLBB for the most general case of \ref{eq:minlp}, which includes nonlinear expressions featuring transcendental functions.  Both \bonmin{} and \minotaur{} provide optimized high-performance C++ implementations of NLBB with a focus on convex \ref{eq:minlp}s.

The core contribution of this work is \juniper{}, a minimalist implementation of NLBB that is designed for rapid exploration of novel NLBB algorithms.  Leveraging the high-level Julia programming language makes it easy to modify \juniper{}'s algorithm and explore extensions, such as branching heuristics, feasibility pumps, and parallelization.  Furthermore, the solver abstraction layer provided by JuMP \cite{DunningHuchetteLubin2017} makes it trivial to change the solvers used internally by \juniper{}'s NLBB algorithm.  Detailed numerical experiments on 300 challenging \ref{eq:minlp}s are conducted to validate \juniper{}'s implementation.  The experiments demonstrate that the initial release of \juniper{} has comparable performance to other established NLBB solvers, such as \bonmin{}, \minotaur{}, and \knitro{}, and that \juniper{} finds high-quality solutions to problems that are challenging for global optimization solvers, such as \couenne{} and \scip{}.  These results illustrate that \juniper{}'s minimalist implementation provides a strong foundation for further exploration of NLBB algorithms for nonconvex \ref{eq:minlp}s.

The rest of the paper is organized as follows. Section \ref{sec:nlbb} provides a brief overview of the NLBB algorithms.  Section \ref{sec:juniper} introduces the \juniper{} NLBB solver.  The experimental validation is conducted in Section \ref{sec:exper}, and Section \ref{sec:conclusion} concludes the paper.

\section{The Core Components of Nonlinear Branch-and-Bound}
\label{sec:nlbb}

To provide context for \juniper's implementation, we begin by reviewing the core components of an NLBB algorithm.  NLBB is a natural extension of the well-known branch-and-bound algorithm for mixed-integer linear programs (MIPs) to \ref{eq:minlp}s.  The algorithm implicitly represents all possible discrete variable assignments in a \ref{eq:minlp} by a decision tree that is exponential in size.  The algorithm then searches through this implicit tree (i.e., {\em branching}), looking for the best assignment of the discrete variables and keeping track of the best feasible solution found thus far, the so-called incumbent solution.  At each node in the tree, the partial assignment of discrete variables is fixed and the remaining discrete variables are relaxed to continuous variables, resulting in a nonlinear program (NLP) that can be solved using an established solver, such as \ipopt{}\cite{ipopt}.   If the solution to this NLP is globally optimal, then it provides a lower bound to the \ref{eq:minlp}'s objective function, $f(x,y)$.  Furthermore, if this NLP {\em bound} is worse than the best solution found thus far, then the NLP relaxation proves that the children of the given node can be ignored.  If this algorithm is run to completion, it will provide the globally optimal solution to the \ref{eq:minlp}.  However, if the \ref{eq:minlp} includes nonconvex constraints, the NLP solver provides only local optimality guarantees, and the NLBB algorithm will be only a heuristic for solving the \ref{eq:minlp}.  The key to designing this kind of NLBB algorithm is to develop generic strategies that find feasible solutions quickly and direct the tree search toward higher-quality solutions.  We now briefly review some of the core approaches to achieve these goals.   

\paragraph{Branching strategy:}
In each node of the search tree, the branching strategy defines the order in which the children (i.e., variable/value pairs) of that node should be explored.  The typical branching strategies are (1) {\em most infeasible}, which branches on the variables that are farthest from an integer value in the NLP relaxation; (2) {\em pseudo cost}, which tracks how each variable affects the objective function during search and then prioritizes variables with the best historical record of improving the objective value \cite{Benichou1971}; (3) {\em strong}, which tests all branching options by brute-force enumeration and then takes the branch with the most promising NLP relaxation \cite{Applegate:1995:FCT:868329}; and (4) {\em reliability}, which uses a threshold parameter to limit strong branching to a specified amount of times for each variable \cite{ACHTERBERG200542}.

\paragraph{Traversal strategy:}
At any point during the tree search there are a number of {\em open} nodes that have branches that remain to be explored.  The traversal strategy determines how the next node will be selected for exploration.  The typical traversal strategies include (1) {\em depth first}, which explores the most recent open node first; and (2) {\em best first}, which explores the open node with the best NLP bound first. The advantage of depth first search is that it only requires a memory overhead that is linear in the number of discrete variables.  In contrast, best first search results in the smallest number of nodes explored but can consume an exponential amount of memory.

\paragraph{Incumbent heuristics:}
In some classes of \ref{eq:minlp}s, finding an initial feasible solution can be incredibly difficult, and the NLBB algorithm can spend a prohibitive amount of time in unfruitful parts of the search tree.  Running dedicated feasiblity heuristics at the root of the search tree is often effective in mitigating this issue.  The most popular such heuristic is the {\em feasibility pump}, which is a fixed-point algorithm that alternates between solving an NLP relaxation of the \ref{eq:minlp} for assigning the continuous variables and solving a MIP projection of the NLP solution for assigning the discrete variables \cite{Fischetti2005,DAmbrosio2012}.

\paragraph{Relaxation restarts:}
In traditional branch-and-bound algorithms, the continuous relaxation is convex and guaranteed to converge to the global optimum or prove that the relaxation is infeasible.  However, in the case of nonconvex \ref{eq:minlp}s, a local NLP solver provides no such guarantees.  Thus, it can be advantageous to restart the NLP solver from a variety of different starting points in the hopes of improving the lower bound or finding a feasible solution \cite{bonmin,knitro}.

\section{The Juniper Solver}
\label{sec:juniper}

The motivation for developing \juniper{} \cite{juniper} is to provide relatively simple and compact implementation of NLBB so that a wide variety of algorithmic modifications can be explored in the pursuit of developing novel heuristics for nonconvex \ref{eq:minlp}s.  To that end, Julia is a natural choice for the implementation for two reasons: (1) Julia provides high-level programming, similar to Matlab and Python, that is preferable for rapid prototyping; and (2) the mathematical programming package JuMP \cite{DunningHuchetteLubin2017} provides an AMPL-like modeling layer, which makes it easy to state \ref{eq:minlp} problems, and a solver abstraction layer, which makes a wide range of NLP and MIP solvers available for use in Juniper.
To demonstrate these properties, Code Block \ref{code:build_solve} provides a simple Julia v0.6 example illustrating the software installation, stating a JuMP v0.18 \ref{eq:minlp} model, and solving it with \juniper{}.  In this example, the NLP solver \ipopt{} is used for solving the continuous relaxation subproblems and the MIP solver \cbc{} is used in the feasibility pump heuristic.

\begin{code}[t]
\caption{Installing and Solving a MINLP with JuMP and \juniper{}}
\label{code:build_solve}
\begin{minted}{julia}
Pkg.add("JuMP"); Pkg.add("Ipopt"); Pkg.add("Cbc"); Pkg.add("Juniper")
using JuMP, Ipopt, Cbc, Juniper

ipopt = IpoptSolver(print_level=0); cbc = CbcSolver()
m = Model(solver=JuniperSolver(ipopt, mip_solver=cbc))

v = [10,20,12,23,42]; w = [12,45,12,22,21]
@variable(m, 0 <= x[1:5] <= 10, Int)

@objective(m, Max, dot(v,x))
@constraint(m, sum(x[i] for i=1:5) <= 6)
@NLconstraint(m, sum(w[i]*x[i]^2 for i=1:5) <= 300)

status = solve(m); getvalue(x) 
\end{minted}
\vspace{-0.4cm}
\end{code}

From Code Block \ref{code:build_solve}, it is clear how \juniper{} can be reconfigured to use different NLP and MIP solvers at runtime.  As is typical for solvers, \juniper{} also features a wide variety of parameters for augmenting the NLBB algorithm.  These include options for selecting the branching strategy, tree traversal strategy, feasibility pump, parallelized tree search, and numerical tolerances, among others.  A complete list of algorithm parameters is available in \juniper{}'s documentation.  After rigorous testing on hundreds of \ref{eq:minlp} problems, the following default settings were identified:  Strong branching is performed at the root node, and pseudo-cost branching is used afterward.  Typically, complete strong branching is conducted; however, if the NLP runtime combined with the number of branches will require more than 100 seconds, the number of branches explored is reduced to meet this time limit.  If the NLP relaxation fails in the root node, it will be restarted up to three times.  Best first search is used for exploring the decision tree, and the runtime of the feasibility pump is limited to 60 seconds.
%

\section{Experimental Evaluation}
\label{sec:exper}

This section conducts a detailed numerical study of \juniper{}'s performance under a variety of configurations and compares its performance to established \ref{eq:minlp} solvers.  Five points of comparison were considered for solving \ref{eq:minlp}s. \bonmin{} v1.8 \cite{bonmin}, \minotaur{} v0.2 \cite{minotaur}, and \knitro{} v10.3 \cite{knitro} were included as alternative NLBB implementations, whereas \couenne{} v0.5 \cite{couenne} and \scip{} v5.0 \cite{scip,scip_v5} were used for a global optimization reference.  All of the open-source solvers utilize \ipopt{} v3.12 \cite{ipopt} compiled with HSL \cite{hsl_lib} for solving NLP subproblems and their respective default LP and MIP solvers.  All of the solvers, except \juniper{}, were accessed through their AMPL NL file interface.  All of the computations were conducted on a cluster of HPE ProLiant XL170r servers featuring two Intel 2.10 GHz 16 Core CPUs and 128 GB of memory.  All solvers were configured with an optimality gap of 0.01\% and a runtime limit of 1 hour.  It is important to note that Julia's JIT takes around 3--10 seconds the first time \juniper{} is run; this time is not reflected in the runtime results.

\paragraph{MINLP problem selection:}

The first step in performing this evaluation is to select an appropriate collection of \ref{eq:minlp} test problems.  We began with 1500 \ref{eq:minlp} problems from MINLPLIB2 \cite{minlplib2}, which are available in Julia via the MINLPLibJuMP package \cite{MINLPLibJuMP}.  Second, all of the problems with no discrete variables or fewer than ten constraints were eliminated, resulting in about 700 problems that focus on the constrained mixed-integer problems that \juniper{} is intended for.  Through a preliminary study, it was observed that more than half of these cases are solved to global optimality or are proven to be infeasible by \scip{} or \couenne{} in less than 60 seconds, suggesting that these are relatively easy cases for state-of-the-art global optimization methods and that they are not of interest to this work. The final collection of test problems consists of 298 \ref{eq:minlp}s that are challenging for both NLBB and global optimization solvers.

\paragraph{Solver comparison:}

\begin{figure}[t]
  \centering
  \vspace{-0.5cm}
  \includegraphics[width=12.0cm]{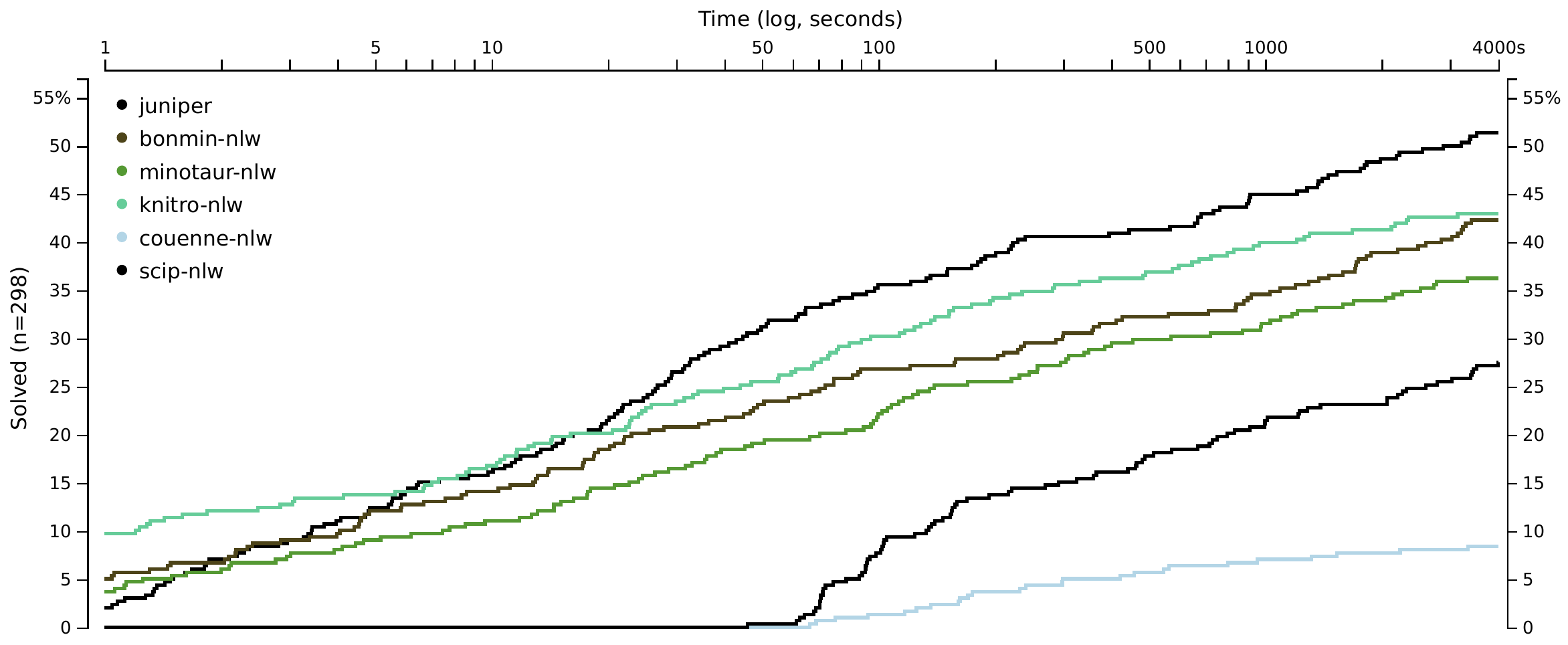}
  \includegraphics[width=6.0cm]{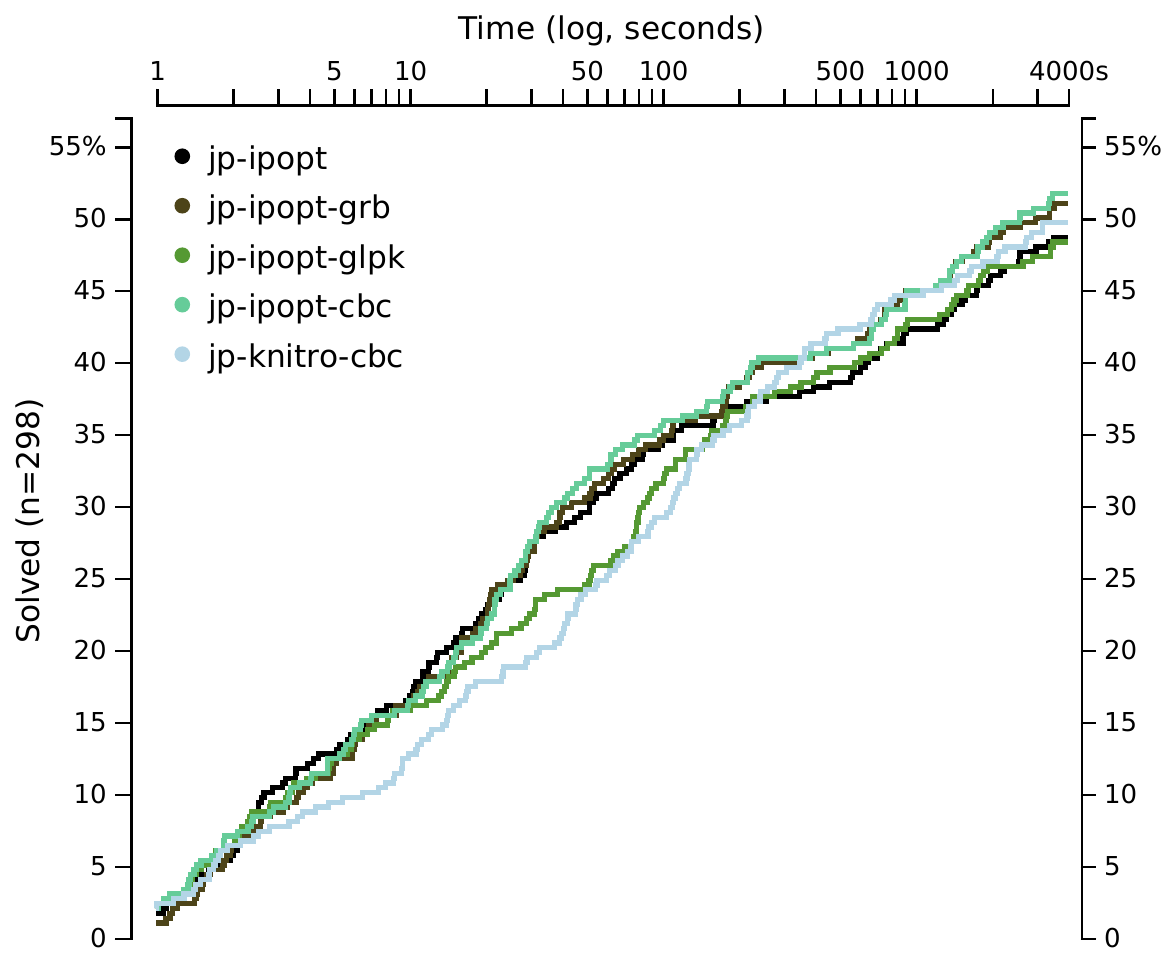}
  \includegraphics[width=6.0cm]{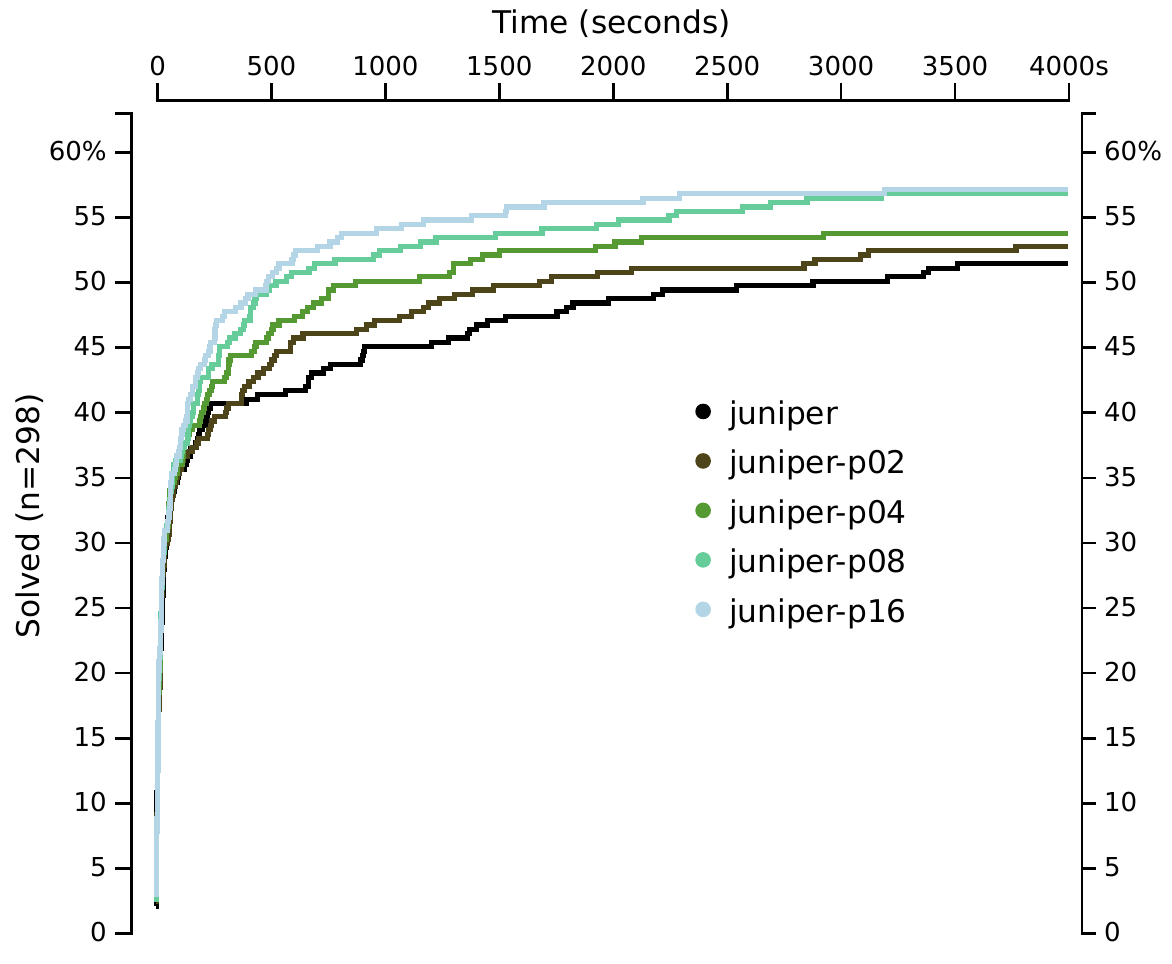}
   \vspace{-0.2cm}
   \caption{Runtime profiles on 298 instances for all different solvers (top) using different MIP solvers (bottom left) and parallelized tree search (bottom right).}
  \label{fig:solvers}
\end{figure}

The first and foremost goal is to demonstrate that \juniper{} has comparable computational performance to \bonmin{} and \minotaur{}. Figure \ref{fig:solvers} (top) provides an overview of the runtime for each solver to complete its tree search procedure.  This figure highlights two key points: (1) \juniper{} is slower for small models that can be solved in less than 30 seconds; however, it consistently solves more models after 30 seconds; and (2) the search completes in no more than 50\% of the cases considered, demonstrating that the selected \ref{eq:minlp} instances present challenging tree search problems for both the NLBB and global optimization solvers.

Table \ref{table:results} provides further details on the performance of each solver, including problem sizes, objective gaps from the best-known solution, and runtime results.  The table begins with summary statistics.  The first row shows the number of feasible solutions found by each solver as well as the number of test cases where the runtime limit was reached. 
The following three rows show the average optimality gaps and runtime for each solver. The first average is for all instances where the specific solver was able to find a feasible solution. The second average is for instances where all six solvers were able to find feasible solutions. The third average is for instances where all four NLBB solvers were able to find feasible solutions.  These summary results indicate two key points: (1) \juniper{} is one of the most robust solvers (only \scip{} had a higher feasible solution count); and (2) for cases where all NLBB algorithms have feasible solutions,  \bonmin{} has the highest solution quality, and \juniper{}, \minotaur{}, and \knitro{} have similar quality, on average.
The remaining rows in the table provide a representative sample of the 298 problems considered.  The first five columns describe each problem by name, number of variables $V$, number of constraints $C$, number of discrete variables $I$, and number of nonlinear constraints $NC$.  The general trends are summarized as follows:
(1) there is a great diversity among which solver is the best on the \ref{eq:minlp} instances considered; and (2) in most cases, the solutions found by the NLBB solvers tie or improve those found by the global solvers; however, there are a few notable cases where global solvers find the best solutions.
Overall, these results indicate that \juniper{} in its default configuration is comparable with the NLBB solvers considered here.


\paragraph{Subsolver selection:}
One of the key features of \juniper{} is that it can use different solvers for the NLP relaxation and for the MIP aspect of the feasibility pump heuristic.
Figure~\ref{fig:solvers} (bottom left) shows a performance profile for a number of subsolver variants of \juniper{}, both with and without a MIP solver (i.e., \glpk{}\cite{glpk}, \cbc{}\cite{cbc}, \gurobi{}\cite{gurobi}), as indicated by \juniper{}-ipopt. \juniper{}-\knitro{}-\cbc{} shows the result of using \knitro{} as the NLP solver instead of \ipopt{}.
The runtime difference between using the feasibility pump and using no heuristic is quite notable in some cases; however, given sufficient time, \juniper{} solves a similar number of cases even without a feasibility pump.  To our surprise, there was little difference in using \cbc{} as the MIP solver compared to using \gurobi{}, suggesting that \cbc{} is a suitable default solver.  We also observed that \glpk{} is not a suitable choice because it was typically unable to terminate in less than 60 seconds, which is the preferred feasibility pump time limit. 


\paragraph{Parallel tree search:}
A key feature of Julia is native and easy-to-use support for parallel processing.  \juniper{} leverages this capability to implement a parallel tree search algorithm.  Figure~\ref{fig:solvers} (bottom right) illustrates the benefits from this simple parallelization of the algorithm, where the first thread orchestrates the computation and all additional worker threads process open nodes in the search tree.  The figure indicates that having two worker threads (instead of using the sequential algorithm) is about 1.7 times faster and having four worker threads is about 3.3 times faster.  The difference between eight and sixteen worker threads is not that notable (both increase speed by around 5.8 times).


\section{Conclusion}
\label{sec:conclusion}

This work has highlighted the potential for leveraging NLBB algorithms as heuristics for solving challenging nonconvex \ref{eq:minlp}s.  To assist in the design of such algorithms, a new Julia-based solver, \juniper{}, is proposed as the base implementation for future exploration in this area.  A detailed experimental study demonstrated that, despite its minimalist implementation, \juniper{} performs comparably to established NLBB solvers on the class of \ref{eq:minlp}s for which it was designed.  We hope that \juniper{} will provide the community with a valuable reference implementation for collaborative open-source research on heuristics for large-scale nonconvex \ref{eq:minlp}s.

\begin{landscape} 
 
\begin{table*}[t] 
\footnotesize 
\caption{Quality and Runtime Results for Various Instances} 
\begin{tabular}{|r|r|r|r|r||r||r|r|r|r|r|r||r|r|r|r|r|r|r|} 
\hline 
 \multicolumn{6}{|c||}{} & juniper    & bon  & minot & knitro & coue        & scip            & juniper          & bon  & minot & knitro  & coue         & scip \\  
    \hline 
    \hline 
\multicolumn{6}{|c||}{Feasible instances / Time limit reached} & 228 & 196 & 193 & 187 & 183 & 257 & 114 & 105 & 134 & 144 & 247 & 215  \\ 
\hline 
\multicolumn{6}{|c||}{} & \multicolumn{6}{c||}{Gap (\%)} &  \multicolumn{6}{c|}{Runtime (seconds)} \\ \hline 
\multicolumn{6}{|c||}{Average solver feasible} & 324.51 & 6.34 & 25.52 & 1365.77 & 395.68 & 1e+04 & 1425 & 1672 & 1789 & 1226 & 3285 & 2853  \\ 
\multicolumn{6}{|c||}{Average all solvers feasible (n=85)} &  21.43 & 0.86 & 3.64 & 24.57 & 16.19 & 2e+04 & 929 & 834 & 1558 & 1002 & 2788 & 2661  \\ 
\multicolumn{6}{|c||}{Average all NLBB solvers feasible (n=113)} &  18.48 & 0.94 & 11.19 & 20.69 & - & - & 918 & 845 & 1420 & 1112 & - & -  \\ 
\hline 
Instance   & $|V|$& $|C|$& $|I|$& $|NC|$ & best obj.  & juniper    & bon  & minot &  knitro & coue        & scip            & juniper          & bon  & minot & knitro  & coue         & scip \\ 
\hline 
                             oil2 &          937 &          927 &           2 &          284 &              -0.73 &  \textbf{0.00} &\textbf{0.00} &  \textbf{0.00} &\textbf{0.00} &              - &             - &            7 &                  3 &         \textbf{2} &                  8 &           - &            - \\ 
                         bchoco08 &          169 &          191 &           9 &           79 &               0.98 &              - &            - &  \textbf{0.00} &         0.71 &              - &             - &            - &                  - &        \textbf{46} &                153 &           - &            - \\ 
                squfl015-080persp &         2416 &         2481 &          15 &         1200 &             402.48 &  \textbf{0.00} &\textbf{0.00} &           4.61 &\textbf{0.00} &           0.93 & \textbf{0.00} &           27 &        \textbf{23} &                259 &        \textbf{23} &         T.L &           94 \\ 
                 transswitch0014p &          139 &          278 &          20 &          121 &            8082.58 &  \textbf{0.00} &\textbf{0.00} &  \textbf{0.00} &\textbf{0.00} &  \textbf{0.00} &             - &            2 &                 14 &         $\bm{< 1}$ &                 32 &         T.L &            - \\ 
                squfl030-100persp &         6031 &         6101 &          30 &         3000 &             363.09 &  \textbf{0.00} &\textbf{0.00} &          43.71 &\textbf{0.00} &         137.24 & \textbf{0.00} &          683 &                919 &                T.L &         \textbf{5} &         T.L &          489 \\ 
                          FLay05M &           63 &           66 &          40 &            5 &              64.50 &  \textbf{0.00} &\textbf{0.00} &  \textbf{0.00} &\textbf{0.00} &  \textbf{0.00} & \textbf{0.00} &         1283 &                561 &               1213 &               1679 &         117 &  \textbf{71} \\ 
                           ndcc13 &          631 &          255 &          42 &           42 &              84.63 &           5.80 &\textbf{0.00} &  \textbf{0.00} &            - &              - &          3.30 &          T.L &               1721 &      \textbf{1172} &                  - &           - &          T.L \\ 
                        CLay0205H &          261 &          366 &          50 &           40 &            8092.50 &  \textbf{0.00} &\textbf{0.00} &         891.18 &         1.61 &              - & \textbf{0.00} &         2548 &               1091 &                T.L &                T.L &           - &  \textbf{71} \\ 
                  procurement1mot &          785 &          750 &          60 &            1 &             291.54 &  \textbf{0.00} &\textbf{0.00} &  \textbf{0.00} &\textbf{0.00} &           8.19 &          1.95 &          445 &       \textbf{204} &                385 &               2345 &         T.L &          T.L \\ 
                      fo9\_ar2\_1 &          181 &          436 &          72 &           18 &              32.62 &              - &            - &         115.90 &            - &          56.21 & \textbf{0.00} &            - &                  - &                T.L &                  - &         T.L &\textbf{3454} \\ 
           crudeoil\_pooling\_ct1 &          311 &          566 &          80 &           37 &          210537.49 &          43.82 &\textbf{0.00} &           3.06 &            - &           2.34 & \textbf{0.00} & \textbf{T.L} &       \textbf{T.L} &       \textbf{T.L} &                  - &\textbf{T.L} & \textbf{T.L} \\ 
               multiplants\_mtg1a &          194 &          257 &          93 &           28 &             391.61 &  \textbf{0.00} &\textbf{0.00} &           3.08 &            - &           2.46 & \textbf{0.00} &          909 &       \textbf{901} &                T.L &                  - &         T.L &          997 \\ 
                multiplants\_mtg2 &          230 &          307 &         112 &           37 &            7099.19 &  \textbf{0.00} &\textbf{0.00} &           0.14 &            - &  \textbf{0.00} &         27.54 &         3390 &      \textbf{1877} &                T.L &                  - &         T.L &          T.L \\ 
           crudeoil\_pooling\_ct3 &          730 &         1199 &         128 &          211 &          287000.00 &          32.69 &            - &              - &            - &          16.82 & \textbf{0.00} & \textbf{T.L} &                  - &                  - &                  - &\textbf{T.L} & \textbf{T.L} \\ 
                    pooling\_epa3 &         1105 &         1718 &         150 &          274 &          -14948.12 &  \textbf{0.00} &            - &              - &         0.17 &              - &             - & \textbf{T.L} &                  - &                  - &       \textbf{T.L} &           - &            - \\ 
                 transswitch0118r &         1240 &         1764 &         179 &         1311 &          129457.85 &           0.45 &            - &              - &\textbf{0.00} &              - &             - & \textbf{181} &                  - &                  - &                491 &           - &            - \\ 
         multiplants\_stg1 &          415 &          262 &         198 &           34 &             355.09 &         24.36 &\textbf{0.00} &           1.03 &             - &              - &             - & \textbf{T.L} & \textbf{T.L} &       \textbf{T.L} &            - &           - &            - \\ 
                              qap &          226 &           31 &         225 &            1 &          389718.00 &           0.50 &         0.38 &           6.79 &\textbf{0.00} &          10.94 &          6.53 &         1404 &                428 &                T.L &       \textbf{137} &         T.L &          T.L \\ 
                          csched2 &          401 &          138 &         308 &            1 &         -166102.00 &              - &\textbf{0.00} &           4.16 &            - &          14.05 &          6.69 &            - &        \textbf{18} &                T.L &                  - &         T.L &          T.L \\ 
                  edgecross22-048 &          463 &         6161 &         462 &            1 &              84.00 &  \textbf{0.00} &\textbf{0.00} &  \textbf{0.00} &\textbf{0.00} &  \textbf{0.00} & \textbf{0.00} &  \textbf{12} &                121 &                 98 &                 23 &         162 &         2201 \\ 
           crudeoil\_pooling\_dt1 &         3713 &         5779 &         916 &          570 &          209585.27 &         368.25 &            - &           4.59 &            - &              - & \textbf{0.00} &          T.L &                  - &                T.L &                  - &           - &\textbf{3429} \\ 
\hline 
\end{tabular}\\ 
\label{table:results} 
\end{table*} 

\end{landscape}

\bibliographystyle{splncs}
\bibliography{references.bib}

\noindent
LA-UR-17-31300

\end{document}